\documentclass{article}
\usepackage{ctex}
\usepackage{amsthm}
\usepackage{amsfonts}
\usepackage{graphicx}
\usepackage{subfigure}
\usepackage{epstopdf}
\usepackage{float}
\usepackage{mathrsfs,makecell}
\usepackage{amstext,graphicx,latexsym,amssymb,amscd,amsfonts,amsmath,slashbox}
\usepackage{float}
\usepackage{color,colortbl,xcolor}
\usepackage{ifthen,calc,lastpage,listings}
\usepackage{amstext, graphicx, latexsym, amssymb, amscd, amsfonts}
\usepackage{booktabs,fancyhdr}
\usepackage{multirow,textcomp,yhmath,exscale,relsize}
\usepackage{ifthen}
\usepackage[left=2.5cm,right=2.5cm,top=2.5cm,bottom=2.5cm]{geometry}
\usepackage[numbers,sort&compress]{natbib}
\usepackage{algorithm}
\usepackage{algorithmic}

\usepackage{amsmath,amsfonts,amssymb,amsthm}
\usepackage{ctex}
\usepackage{graphicx}
\usepackage{float}
\usepackage{geometry}
\usepackage{bm}
\newcolumntype{L}[1]{>{\flushleft}p{#1}}
\newcolumntype{C}[1]{>{\centering}p{#1}}

\begin{document}
 \def\QED{\hfill $\Box$\smallskip}

 \begin{center}
 {\Large\bf The  Quasi-Newton Method for the Composite Multiobjective Optimization Problems}
 \noindent \footnote{This research  was partially supported by   the National Natural Science Foundation of China (12271071, 11991024),   the Team Project of Innovation Leading Talent in Chongqing (CQYC20210309536), ¡°Contract System¡± Project  of Chongqing Talent Plan (cstc2022ycjh-bgzxm0147), the Chongqing University Innovation Research Group Project (CXQT20014) and  the Basic and Advanced Research Project of
Chongqing ( cstc2021jcyj-msxmX0300).}\\

 \vspace{12pt}
 \noindent
Jian-Wen Peng\footnote{ Corresponding author. School of Mathematical Sciences, Chongqing Normal University, Chongqing 401331, China. e-mail: jwpeng168@hotmail.com}
  and
  Jen-Chih Yao\footnote{ Department of Applied
 Mathematics, National Sun Yat-sen University, Kaohsiung, Taiwan 804.
 e-mail: yaojc@math.nsysu.edu.tw}
\vspace{30pt}

 \end{center}

 \clearpage

 \noindent {\bf Abstract}  In this paper, we introduce several new    quasi-Newton methods for the composite multiobjective optimization problems (in short, CMOP) with Armijo line search. These multiobjective versions of   quasi-Newton methods  include   BFGS  quasi-Newnon method,  self-scaling BFGS    quasi-Newnon method, and Huang BFGS    quasi-Newnon method.   Under some suitable conditions,  we show that each accumulation point of the sequence generated by these  algorithms, if exists, is  both a Pareto stationary point and  a Pareto optimal point of (CMOP).

 \vspace{12pt}

 \vspace{12pt}

\clearpage

\vskip 0.4in
\renewcommand{\theequation}{\arabic{equation}}
\section{Introduction} \vskip 0.2in

 \noindent In this work, let us consider the    following  composite   multiobjective optimization problems (in short, CMOP):
\begin{equation}
\begin{aligned} \label{P}
\quad &\min \quad F(x)  \\
&\begin{array}{r@{\quad}r}
s.t. & x \in { \mathbb{R}^n}, \\
\end{array}
\end{aligned}
\end{equation}
where $F:{\mathbb{R}^n} \to \mathbb{R} $ is a vector-valued function with $F: = {({f_1},..,{f_{\rm{m}}})^T}$ and $T$ denotes transpose.  We assume that each ${f_i}:{\mathbb{R}^n} \to \mathbb{R} $ is defined by
\begin{align}
{f_i}(x): = {g_i}(x) + {h_i}(x), ~~~ i = 1,...,m,
\end{align}
where ${g_i}{\rm{:}}{\mathbb{R}^n} \to \mathbb{R}$ is a twice continuously differentiable strongly convex function, ${h_i}:{\mathbb{R}^n} \to \mathbb{R} $ is  convex  but not necessarily differentiable.
It is worthy noting that if $h_i(x) \equiv 0$ for all $x\in {\mathbb{R}^n}$ and $i=1, 2, ...,m$, then  the above  (CMOP) (i.e.,  (1) with (2) ) reduces to the    multiobjective optimization problems studied in \cite{FG, FS, BLM, BIS, DDFL, FGS, P}.

Scalarization approach is one of the most popular strategies  to solve the multiobjective optimization problem, which transforms the multiobjective optimization problem into one or several parameterized single-objective  problem (see \cite{EG, EM,NYY} and the references therein). In general, the converted problem
and the primal multiobjective optimization problem enjoy same optimal solutions under certain conditions. Nevertheless, Fliege el al. \cite{FGS}
pointed out that the parameters in scalarization method are not known in advance and the selection of
parameters may result in unbounded scalar problems even if the original multiobjective optimization problem has solutions. In order to
cope with these limitations, the descent methods for multiobjective optimization problems has attracted wide attention in the optimization field (see \cite{FG}).  There are many kinds of descent methods for multiobjective optimization problems, for examples,   the steepest descent method   \cite{FS},  the projection gradient method \cite{BLM}, the proximal point method  \cite{BIS}, the  Newton's method \cite{FGS}, quasi-Newton method \cite{ P,  QGC, MBP}, and  the subgradient method \cite{DDFL}.
Recent years, some authors introduced and researched the descent methods for the  composite   multiobjective optimization problem (1.1), for examples,     the proximal gradient method \cite{TFY} ,    the  Newton's method \cite{Ansary},  proximal Newton methods  \cite{RP},  proximal  quasi-Newton methods \cite{PRY}.   Inspired by above   ideas,     the main purpose of this paper is to  introduce  a new     quasi-Newton methods for    (CMOP)     in the   case   of Armijo line search.

The main contents of this paper are as follows: In Section 2, we give some notations and some concepts about Pareto optimality and Pareto stationarity. In Section 3,   we propose BFGS   quasi-Newnon method, self-scaling BFGS   quasi-Newnon method    and  Huang BFGS   quasi-Newnon method      for the (CMOP) with  the Armijo line search.   We prove the global convergence of the proposed algorithms in Section 4.  In Section 5,   some numerical experiments are also carried out to verify the effectiveness of the proposed   quasi-Newton methods and to show that  multiobjective version of the Huang BFGS proximal  quasi-Newnon method    is the   most   effective     among  the   gradient method,   Newton method and the proposed  quasi-Newton methods for (CMOP) with  an Armijo line search.

\vskip 0.4in

\section{ Preliminaries}
\vskip 0.2in

Throughout this paper,  $\mathbb{N}$ denotes the set of positive integers, $\mathbb{R}^n$ denotes  the $n$-dimentional Euclidean space. The Euclidean norm in ${\mathbb{R}^n}$ will be denoted by $\left\|  \cdot  \right\|$. For two vectors $u$ and $v$ in $\mathbb{R}^m$, by    $u \le v$, we mean  that ${u_i} \le {v_i} $ for all $i = 1,...,m.$  by    $ u < v $, we mean  that $ {u_i} < {v_i}$ for all $i = 1,...,m.$

  ~~We recall   that  the function $g:{\mathbb{R}^n} \to \mathbb{R}$ is   strongly convex or uniformly convex if  there exists  $\eta>0$ such that
$$ g(\lambda x+ (1-\lambda)y)\leq \lambda g(x) +(1-\lambda)g(y) -\frac{1}{2} \lambda (1-\lambda)\eta \left\| {x - y} \right\|^2,$$
 for any $ x,y \in {\mathbb{R}^n}$ and $\lambda\in [0,1]$  (see \cite{NW}).
 It is well known that if $g:{\mathbb{R}^n} \to  \mathbb{R}$ is  twice continuously differentiable, then   $g$  is strongly convex if and only if    there exists  $\eta>0$ such that
       $   {\nabla ^2}g(x) \succeq \eta I$, for any $x \in {\mathbb{R}^n}$,
where ${\nabla ^2}g(x)$ denote the Hessian matrix of $f$ at $x$. Hence£¬if ${g}$ is strongly convex, then the Hessian matrix ${\nabla ^2}{g}(x)$ is positive definite for any $x \in {\mathbb{R}^n}$. It is clearly that the  strong convexity of $g$ implies  its strict convexity and usual convexity.

Let $f{\rm{:}}{\mathbb{R}^n} \to \mathbb{R} \cup \{  + \infty \} $, and let $x\in dom(f):=\{x\in \mathbb{R}^n: f(x)< +\infty\}$. Then the directional derivative of $f$ at $x$ in the direction $d\in {\mathbb{R}^n}$ is defined to be the limit

\[f'(x;d) = \mathop {lim}\limits_{\alpha  \to {0^ + }} \frac{{f(x + \alpha d) - f(x)}}{\alpha },\]
if it exists (see \cite{R}). Clearly, if $f$ is differentiable at $x$, then $f'(x;d) = \nabla f{(x)^T}d$.

\textbf{Definition  2.1} \cite{FG, FS}~~Recall that ${x^*} \in {{\mathbb{R}^n}}$ is an efficient  point or  a Pareto optimum of (CMOP), if there is no $x \in {{\mathbb{R}^n}}$ such that $F(x) \le F({x^*})$ and $F(x) \ne F({x^*})$. The set of all Pareto optimal values is called Pareto frontier. Likewise, ${x^*} \in {{\mathbb{R}^n}}$ is a weakly  efficient  point or  a weakly Pareto optimum of  (CMOP), if there is no $x \in  {\mathbb{R}^n}$ such that $F(x) < F({x^*})$.

It's well known that the   efficient  point of (CMOP)  is also  a weakly   efficient  point   of (CMOP), and the converse is not  true in general.

\textbf{Definition 2.2} \cite{TFY}~ We say that $\bar x \in  {\mathbb{R}^n}$ is Pareto stationary (or critical) of (CMOP), if and only if
\[\mathop {\max }\limits_{i = 1,...,m} {f'_i}(\bar x;d) \ge 0~ ~{\rm for ~all ~}~d \in  {\mathbb{R}^n} .\]

It is worthy to noting that Definition 2.2 generalizes the corresponding ones in \cite{FS}.

 We recall the following important result about the relationship of the weakly   efficient  point (efficient  point)  and the  Pareto stationary point of (CMOP):\\

\textbf{Lemma  2.1 {\cite{FGS}}}   ~~(1) If $x \in {R^n}$ is a weakly   efficient  point of (CMOP)£¬then $x$ is Pareto stationary.

                          (2) Let every component ${f_i}$ of $F$ be convex. If $x \in {R^n}$ is a Pareto stationary point of (CMOP), then $x$ is a weakly   efficient  point of (CMOP).

                          (3) Let every component ${f_i}$ of $F$ be strictly convex. If $x \in {R^n}$ is a Pareto stationary point of (CMOP), then $x$ is an   efficient  point of (CMOP).

\vskip 0.4in

\section{ Quasi-Newton methods for (CMOP)}
\vskip 0.2in

In this section, we consider a new    quasi-Newton method  for (CMOP)   with an Armijo line search. Prior to that, we need the following result, which was shown in \cite{Ansary}.\\

\textbf{Lemma 3.1}    Suppose that $0\in Co{ \mathop{ \cup} \limits_{j\in \{{i = 1,...,m}\}} \partial f_j(x^*)}$ for some $x^*\in {\mathbb{R}^n}$, then $x^*$ is a critical point of the (CMOP).\\

 For $i=1,2,...,m,$ let  $\nabla {g_i}(x)$ denote  the gradient of ${g_i}$ at $x$, ${B_i}(x)$ be an    an approximation of ${\nabla ^2}{g_i}(x)$, which   satisfies  the following assumption:

{\bf Condition 3.1} For    any fixed $x\in \mathbb{R}^n$ and for each $i\in \{1, 2,...,m\}$, there exists a constant $\sigma_i$ such that   $$z^T {B_i}(x) z\geq \sigma_i, \forall z\in \mathbb{R}^n.$$

{\bf Remark 3.1}  (i)  Condition 3.1 implies that for any $x\in \mathbb{R}^n$ and for each $i\in \{1, 2,...,m\}$, there exists a constant $\sigma$ such that   $$z^T {B_i}(x) z\geq \sigma, \forall z\in \mathbb{R}^n,$$
where $\sigma= \mathop {\min }\limits_{i = 1,...,m} \sigma_i $.\\

(ii) If ${B_i}(x)=\nabla^2 g_i (x)$, then  Condition 3.1 holds true automatically because the strong convexity of $g_i$.

\vskip 0.4in

 Now, We define the function
 $ \theta(.,.):{\mathbb{R}^n} \times{\mathbb{R}^n}\to \mathbb{R} $ by

\begin{align}
{ \theta(x,d): = \mathop {\max }\limits_{i = 1,...,m} \theta_i(x,d)},
\end{align}
where $\theta_i(x,d)=  \nabla {g_i}{(x)^T}d + \frac{1}{2}{d^T}{B_i}(x)d + {h_i}(x + d) - {h_i}(x). $

It is easy to see that for any fixed $x$, $ \theta(x,d)$ is condinuous on $d$ since $\theta_i(x,d)$ is  condinuous on $d$ for each $i\in \{1, 2,...,m\}$.    It follows from Condition 3.1 and the convexity of  ${h_i}$, we know that  for any fixed $x$ and for each $i\in \{1, 2,...,m\}$,  $\theta_i(x,d)$ is $\sigma_i$-strongly convex in $d$.
Therefore, for any fixed $x$, $\theta(x,d)$  is $\sigma$-strongly convex in $d$ and $ \theta(x,\bm{0})=0$, where $\sigma= \mathop {\min }\limits_{i = 1,...,m} \sigma_i $.  We would like to define the quasi-Newton direction at an iteration $k$  as ${d_{QN}^k} = {d}({x^k})$,   where
\begin{align}
{d }(x): = \mathop {argmin}\limits_{d \in {\mathbb{R}^n}} {\theta(x,d)}.
\end{align}

For any fixed $x$, we can rewrite problem (4) as  the following subproblem $\Phi(x)$, which is to find a suitable descent direction of the (CMOP):
\begin{align}
\Phi(x): \mathop {\min_{d\in {\mathbb{R}^n} }} \theta(x,d).\end{align}

\textbf{Remark 3.2}~~(i) For fixed $x\in  {\mathbb{R}^n}$, it follows from the strongly convexity of ${\theta(x,d)}$ in $d$ that (5) (i.e., $\Phi(x)$) has  a unique solution ${d}(x)$. Moreover,      the optimal value of  $\Phi(x)$ is denoted by  $\alpha(x)$, i.e.,
\begin{align} \alpha(x)= \mathop {\min_{d\in {\mathbb{R}^n} }} \theta(x,d)= \theta(x,{d}(x)).\end{align}

                   (ii) For   fixed point $x\in {\mathbb{R}^n}$, because of  $\theta(x,\bm{0}) = 0$, we have ${\theta(x, d(x))} \le 0.$\\

  \textbf{Remark 3.3} (i) ${\theta(x,d)}$ is  denoted by ${\theta _{x}}(d)$ in  \cite{PRY} and  $ \theta(x,d)$ is exactly the ${\varphi _{\omega ,x}}(d)$  in  \cite{PRY}  with $\omega=0$.

                   (ii) If  for $i=1,2,...,m$, ${B_i}(x)$ are replaced by   ${\nabla ^2}{g_i}(x)$, then $ \theta(x,d)$ reduces to $Q(x, d)$ in \cite{Ansary}  and
                    ${d_{QN}^k} = {d}({x^k})$ defined by (5) reduces to exactly the Newton direction at an iteration $k$ in \cite{Ansary} which is defined as    ${d_{N}^k}=:  \mathop {argmin} \limits_{d \in {\mathbb{R}^n}} Q(x^k, d)$.

We recall an important property of ${\theta(x, .)}$ as follows:

\textbf{Lemma 3.2 ( see \cite{PRY}})~~For   fixed point $x\in {\mathbb{R}^n}$ and for all $d \in {\mathbb{R}^n}$, the following equality holds:
$$ \theta'((x, \bm{0}); d) = \mathop {\max}\limits_{i = 1,...,m} {f'}_i(x;d).$$
$~~$\\

Denote $I(x, d)= \{ j\in \{ 1, 2, ...,m\} : \theta(x,d)=\theta_j(x,d)\}.$
It follows that ${d }(x)= \mathop {argmin}\limits_{d \in {\mathbb{R}^n}} {\theta(x,d)}$ and $\alpha(x)= \theta(x,{d}(x))\leq  \theta(x, \bm{0}) = 0$  that for every $x\in {\mathbb{R}^n}$, $d(x)$ is a solution of $\Phi(x)$). Therefore, $$ \bm{0}\in \partial_d \theta(x,{d}(x)).$$

It follows from Corollary 3.5 in \cite{Bagirov} that there exist  $w\in {\mathbb{R}_+^{|I(x, d(x)) |}}$, $\xi_j\in \partial_d h_j(x+d(x)), j\in I(x, d(x))$ such that the following conditions hold:

\begin{align}
\sum_{j\in I(x,d(x)) } w_j=1.\end{align}

\begin{align}
\sum_{j\in I(x,d(x)) } w_j ( \nabla g_j(x) +B_j(x) d(x) +\xi_j )   =0.\end{align}

Let $\Xi_m:=\{1, 2, ...,m\}$ and substituting $w_j =0$ and $\xi_j\in \partial_d h_j(x+d(x))$ for all $j\notin I(x, d(x)$, we can obtain

\begin{align}
 \sum_{j=1}^m    w_j=1.\end{align}

\begin{align}
\sum_{j=1}^m w_j ( \nabla g_j(x) +B_j(x) d(x) +\xi_j )   =0.\end{align}

\begin{align}
w_j \geq 0,~ w_j(   \nabla g_j(x)^T d(x) +\frac{1}{2}  d(x)^T B_j(x) d(x) +  h_j(x+d(x))-h_j(x)- \alpha(x) )   =0, j\in \Xi_m.\end{align}

\begin{align}
     \nabla g_j(x)^T d(x) +\frac{1}{2}  d(x)^T B_j(x) d(x) +  h_j(x+d(x))-h_j(x) \leq  \alpha(x), j\in \Xi_m.\end{align}

 Therefore, we have the following result:\\

  {\bf Lemma 3.3} If $d(x)$ is a solution of $\Phi(x)$ and $ \alpha(x)=  {\theta(x, d(x))}$, then

   (i) there exist  $w\in {\mathbb{R}_+^{|I(x, d(x)) |}}$, $\xi_j\in \partial_d h_j(x+d(x)), j\in I(x, d(x))$ such that
   such that $(d(x), \alpha(x), w)$ satisfies (7)-(8).\\

(ii)  there exists $w\in {\mathbb{R}_+^{m}}$ such that $(d(x), \alpha(x), w)$ satisfies (9)-(12).\\

The following lemma characterizes the Pareto stationarity of (CMOP) in terms of ${d }( \cdot )$.\\

\textbf{Lemma 3.4}~~Suppose that Condition 3.1 holds true,  ${d}(x)$ and ${\alpha}(x)$ are defined in (4) and (6), respectively. Then, the following statements hold true:

                (1) If $x$ is a Pareto stationary point of (CMOP), then ${d}(x) = 0$ and ${\alpha}(x) = 0.$ Conversely, if ${d }(x) = 0$ and ${\alpha}(x) = 0$, then $x$ is a Pareto stationary point of (CMOP).

                (2) If $x$ is not a Pareto stationary point of (CMOP), then ${d}(x) \ne 0$ and $\alpha(x) < 0$. Conversely, if ${d }(x) \ne 0$ and ${\beta  }(x) < 0$, then $x$ is not a Pareto stationary point of (CMOP).

\begin{proof}[\quad\quad\textbf{Proof}]
(1)  Let $x$ be Pareto stationary of (CMOP). We will prove that  ${d }(x) = 0$ and $\alpha(x) = 0.$  On the contrary,  we assume, that ${d }(x) \ne 0$ or ${\alpha}(x) < 0.$ It follows from   Remark 3.2  that ${d }(x) \ne 0$ if and only if ${\alpha}(x) < 0$, which means that ${d }(x) \ne 0$ and ${\alpha}(x) < 0.$ It follows from (12) in Lemma 3.3(ii) and the positiveness  of $B_j$ that

\begin{align}
     \nabla g_j(x)^T d(x)  +  h_j(x+d(x))-h_j(x)\leq \alpha(x)-\frac{1}{2}  d(x)^T B_j(x) d(x)<0,  j\in \Xi_m.\end{align}

It follows from the convexity of $h_j$ that for any $\lambda \in (0, 1)$ that
\begin{align}
       h_j(x+\lambda  d(x))-h_j(x)\leq  & \lambda (h_j(x+   d(x)) +(1- \lambda) h_j(x)- h_j(x) \notag \\
       = & \lambda [h_j(x+   d(x)) - h_j(x)],  j\in \Xi_m.\end{align}

 From  (13) and (14), we obtain that

$$
    \lambda  \nabla g_j(x)^T d(x) +  h_j(x+\lambda d(x))-h_j(x)
     <  0,  j\in \Xi_m,$$
which implies that

$$
 \frac{1}{\lambda} [ \lambda \nabla g_j(x)^T d(x) +  h_j(x+\lambda d(x))-h_j(x)]
     <0,  j\in \Xi_m. $$

Taking limit $\lambda\rightarrow 0^+$ in the above inequality we have $  {f'_j}(x;d(x))<0$ for each $j= 1, 2, ...,m$. Hence, $x$ is not a critical point of (CMOP), which is a contradiction.

We now  prove the converse. We assume that $d(x)=\bm{0}$ and ${\alpha}(x) = 0$ hold true.  It follows from Lemma 3.3(i) and $ I(x, \bm{0})= \{1,2,...,m\}=\Xi_m$ that  there exist  $w\in \mathbb{R}_+^{m}$, $\xi_j\in \partial_d h_j(x+d(x)), j\in \Xi_m$ such that
   such that $(d(x), \alpha(x), w)$ satisfies
   $$\sum_{j=1}^m  w_j=1, {\rm~~and~~  } \sum_{j=1}^m  w_j ( \nabla g_j(x) +\xi_j )   =0,$$
which implies that $$\bm{0}\in Co{ \mathop{\cup} \limits_{ j\in \Xi_m}\partial f_j(x)}.$$
Therefore, it follows from Lemma 3.1 that $x$ is a critical point of (CMOP). \\

(2) This statement is equivalent to statement (1).
\end{proof}

\textbf{Lemma 3.5}  ~~Let $d (x)$ and ${\alpha}(x)$ be defined in (4) and (6), respectively. Then,  the mappings ${d}( \cdot )$ and ${\alpha}( \cdot )$ are continuous.

\begin{proof}[\quad\quad\textbf{Proof}] Clearly, the following function
\[{ \theta(x,d): = \mathop {\max }\limits_{i = 1,...,m} \theta_i(x,d)}=\mathop {\max }\limits_{i = 1,...,m} \{ \nabla {g_i}{(x)^T}d + \frac{1}{2}{d^T}{B_i}(x)d + {h_i}(x + d) - {h_i}(x)\}  \]
is continuous with respect to $x$ and $d$. Therefore, it follows from \cite[Proposition 23]{AE} that the optimal value function $\alpha( \cdot )$ is also continuous. Furthermore, since the optimal set mapping ${d}( \cdot )$ is unique,  it follows from   \cite[Corollary 8.1]{H} that ${d}( \cdot )$ is continuous.
\end{proof}

 From Lemma 3.5, it is easy to see that the following result holds true.\\

{\bf Corollary 3.1}  Let  $x_k\in  \mathbb{R}^n$ and  $d(x_k)$  be the solution of $\Phi(x_k)$.

 (i) suppose that $\{x_k\}$ converges to $x^*$ and $d(x^k)$ converges to  $d^*$, then $d^*=d(x^*)$.

 (ii)suppose that $\{x_k\}$ converges to $x^*$ and $\alpha(x^k)$ converges to  $\alpha^*$, then $\alpha^* =\alpha(x^*)$.\\

\textbf{Remark 3.4}~(i)~From Lemma 3.5, we know that  $x$ is a Pareto stationary point of (CMOP) if and only if ${d}(x) = 0$, which extends and improves Lemma 3.2 in \cite{Ansary} since the $\nabla^2 g_j(x)$ has been replaced by a   positive definite matrix $B_j(x)$ for any $j\in \Xi_m$.

(ii) Corollary 3.1 extends and improves Theorem  3.1 in \cite{Ansary}.

\textbf{Theorem 3.1} Suppose that Condition 3.1 holds true.
 
 (i) If   $d(x)$ and $\alpha(x)$ are the   solution and  the optimal value  of $\Phi(x)$, respectively, then
\begin{align} \alpha(x) \leq -\frac{\sigma}{2} \|d(x)\|^2,\end{align}
where $\sigma= \mathop {\min }\limits_{i = 1,...,m} \sigma_i$. 

(ii)   if $\rho\in (0,1)$ and  $x$ is a non critical point of (CMOP), then the following inequality    holds  true for every $\lambda > 0$ sufficiently small,
\begin{align}
     f_j(x+\lambda d(x)) \leq f_j(x) + \rho \lambda  \alpha(x)\end{align}
holds for any $j\in \Xi_m$.

{\bf Proof.} (i) Suppose that    $d(x)$ is a solution of $\Phi(x)$ and $ \alpha(x)= \theta(x, d(x))$. Then there exists $w\in {\mathbb{R}_+^{m}}$ such that  $(d(x), \alpha(x), w)$ satisfies (9)-(12). Since $h_j$ is convex and $\xi_j\in \partial_d h_j(x+d(x))$,
\begin{align} h_j(x+d(x))- h_j(x) \leq \xi_j^T d(x).\end{align}
Multiplying both sides of (10) by $d(x)$,
$$
\sum_{j=1}^m  w_j [ \nabla g_j(x)^T d(x) + d(x)^T B_j(x) d(x) +\xi_j^T d(x) ]   =0.$$
 It follows from  the above equality and (17) that
\begin{align}
\sum_{j=1}^m   w_j [ \nabla g_j(x)^T d(x) + d(x)^T B_j(x) d(x) + h_j(x+d(x))- h_j(x)] \leq 0.\end{align}
Taking sum over $j\in \{1, 2, ...,m\}$ in (11) and using (9),
$$
\sum_{j=1}^m   w_j [ \nabla g_j(x)^T d(x) + d(x)^T B_j(x) d(x) + h_j(x+d(x))- h_j(x) ]
=\sum_{j=1}^m   w_j \frac{1}{2}  d(x)^T B_j(x) d(x) + \alpha(x).$$

It follows from  (18) that
\begin{align}
\alpha(x)  \leq  - \sum_{j=1}^m   w_j \frac{1}{2}  d(x)^T B_j(x) d(x).\end{align}

From (19), (9),  Condition 3.1 and the definition of $\sigma$, we obtain
$$ \alpha(x) \leq -\frac{\sigma}{2} \|d(x)\|^2.$$

(ii) Suppose that $x$ is non critical. Then from Lemma 3.4, $d(x)\neq 0$. It follows from (15) that $\alpha(x)<0$.

From Condition 3.1 and the convexity of $h_i$, we know that  for $i\in \Xi_m$ and for  any $\lambda\in [0, 1]$,

\begin{align*}
{f_i}({x} + {\lambda }{d(x)}) -{f_i}({x})
&={g_i}({x} + {\lambda }{d(x)}) -{g_i}({x}) +{h_i}({x} + {\lambda }{d(x)}) -{h_i}({x})\\
&=\lambda (\nabla g_i(x)^Td(x)  +{h_i}({x} + {\lambda }{d(x)}) -{h_i}({x}))+ o(\lambda)\\
&<     \lambda (\nabla g_i(x)^Td(x) +\frac{1}{2} d(x)^T B_i(x) d(x) +{h_i}({x} + {\lambda }{d(x)}) -{h_i}({x}) ) + o(\lambda)\\
&\le  \lambda \alpha(x) + o(\lambda).
\end{align*}

 Then for $i\in \Xi_m$,
$$
{f_i}({x} + {\lambda }{d(x)}) -{f_i}({x})-\lambda \rho \alpha(x) \leq  \lambda (1-\rho) \alpha(x) + o(\lambda).$$
Since $\rho\in (0,1)$ and $ \alpha(x)<0$, the right hand side term in the above inequality becomes non positive for every $\lambda>0$ sufficiently small, which implies that (16) holds true for every $\lambda>0$ sufficiently small.

$~~~~~$\\
{\bf Rewmark 3.5.} If  $B_j(x)\equiv \nabla^2 g_j(x)$  for any $j\in \Xi_m$, then by  Theorem 3.1 we recover    Theorem 3.2 in \cite{Ansary}.\\

To solve  for (CMOP), now we present our new  BFGS quasi-Newton methods with line searches. We compute the step length ${\lambda _k} > 0$ by   an Armijo rule.
Let $\rho   \in (0,1)$ be a   constant. If for each $i \in \Xi_m$, the following inequality holds true
\begin{align}
{f_i}({x^k} + {\lambda _k}{d^k}) \le {f_i}({x^k}) + {\lambda _k}\rho \alpha(x^k),
\end{align}
then the ${\lambda _k}$ is  accepted. Otherwise, we begin with ${\lambda _k} = 1$ and if there exists $i\in \{1, 2, ...,m\}$ such that   the   inequality  (20) is not satisfied, we update
\[{\lambda _k}: = \zeta {\lambda _k},\]
where $\zeta  \in (0,1)$.

The following result illustrates that  ${d^k}$ produced by the Armijo rule  procedure is a descent direction of (CMOP) at a nonstationary points  ${x^k}$.\\

\textbf{Lemma 3.6}~Suppose that Condition 3.1 holds true. ~Let $\rho  \in (0,1)$, ${d^k}:= d(x^k)$ and $\alpha(x^k)$ be the solution and the optimal value of $\Psi(x^k)$, respectively. If ${x^k}$ is not Pareto stationary, then there exists some ${\bar \lambda _k} > 0$ such that for each $i = 1,...,m$  and for any $\lambda  \in (0,{\bar \lambda _k}]$, the following inequality holds true
\[{f_i}({x^k} + \lambda {d^k}) \le {f_i}({x^k}) + \lambda \rho \alpha(x^k).\]

\begin{proof}[\quad\quad\textbf{Proof}]
Let $\lambda  \in (0,1].$ It follows from the  convexity of  ${h_i}$ that  for each $i = 1,...,m,$ we have
\begin{align}
{h_i}({x^k} + \lambda {d^k}) - {h_i}({x^k}) = & {h_i}((1 - \lambda ){x^k} + \lambda ({x^k} + {d^k})) - {h_i}({x^k}) \notag\\
\le & (1 - \lambda ){h_i}({x^k}) + \lambda {h_i}({x^k} + {d^k}) - {h_i}({x^k}) \notag \\
= & \lambda [{h_i}({x^k} + {d^k}) - {h_i}({x^k})]. \notag
\end{align}

From Condition 3.1 and  the first-order Taylor expansion of ${g_i}$, we have that for  each $i\in \Xi_m,$
\begin{align}
&{g_i}({x^k} + \lambda {d^k}) + {h_i}({x^k} + \lambda {d^k}) \notag \\
\le & {g_i}({x^k}) + \lambda \nabla {g_i}{({x^k})^T}{d^k} + \frac{1}{2}{(\lambda {d^k})^T}{B_i}({x^k})(\lambda {d^k}) + {h_i}({x^k}) + \lambda ({h_i}({x^k} + {d^k}) - {h_i}({x^k})) + o({\lambda}) \notag \\
= & {g_i}({x^k}) + {h_i}({x^k}) + \lambda [\nabla {g_i}{({x^k})^T}{d^k} + \frac{\lambda }{2}{({d^k})^T}{B_i}({x^k})({d^k}) + {h_i}({x^k} + {d^k}) - {h_i}({x^k})] + o({\lambda}) \notag \\
\le & {g_i}({x^k}) + {h_i}({x^k}) + \lambda [\nabla {g_i}{({x^k})^T}{d^k} + \frac{1}{2}{({d^k})^T}{B_i}({x^k})({d^k}) + {h_i}({x^k} + {d^k}) - {h_i}({x^k})] + o({\lambda}) \notag \\
\le & {g_i}({x^k}) + {h_i}({x^k}) + \lambda \alpha(x^k) + o({\lambda}) \notag \\
= & {g_i}({x^k}) + {h_i}({x^k}) + \lambda \rho \alpha(x^k) + \lambda \left[ {(1 -\rho)\alpha(x^k) + \frac{{o({\lambda})}}{\lambda }} \right], \notag
\end{align}
where ${B_i}({x^k})$ is some approximation of ${\nabla ^2}{g_i}({x^k}),i\in \Xi_m,$ the second inequality follows from the positive definiteness of ${B_i}({x^k})$ and $\lambda  \in (0,1]$, and the third one comes from the definition of $\alpha(x)$ in Remark 3.2. Since ${x^k}$ is not Pareto stationary, we have $\alpha(x^k) < 0$ from Lemma 3.5. It follows from $\rho \in (0,1)$ that there exists some ${\bar \lambda _k} > 0$ such that for  each $i \in \Xi_m,$
\[{g_i}({x^k} + \lambda {d^k}) + {h_i}({x^k} + \lambda {d^k}) \le {g_i}({x^k}) + {h_i}({x^k}) + \lambda \rho \alpha(x^k),\]
 for any  $\lambda  \in (0,{\bar \lambda _k}]$. \end{proof}

To simplify the notation we will use $B_i^k$ to denote ${B_i}({x^k})$  for all $i = 1,...,m$ and $k=0, 1,2,...$.

Now,  we would like to state our new   quasi-Newton methods with line searches for (CMOP) as follows:

\textbf{Algorithm 3.1}

Step 1 Choose $\omega  > 0$, $\rho  \in (0,1)$, $\zeta  \in (0,1)$, ${x^0} \in {\mathbb{R}^n}$, symmetric positive definite matrix $B_i^0 \in {\mathbb{R}^{n \times n}},i = 1,...,m$ and set $k: = 0$;

Step 2 Compute $ d(x^k)$ and $\alpha(x^k)$ by solving subproblem (4) with $x = {x^k}$, let ${d^k}:= d(x^k)$;

Step 3 If ${d^k} = 0$, then stop. Otherwise, proceed to the next step;

Step 4 Compute the step length ${\lambda _k} \in (0,1]$ as the maximum of the following set
\[{\Lambda _k}: = \{ {\lambda } = {\zeta ^j}|j \in \mathbb{N},{f_i}({x^k} + {\lambda }{d^k}) \le {f_i}({x^k}) + {\lambda }\rho \alpha(x^k),i = 1,...,m\} ;\]

Step 5 Set ${x^{k + 1}} = {x^k} + {\lambda _k}{d^k}$, update $\{B_i^{k }\}$  by   following  BFGS update  formula for each $i=1, 2, ...,m$
\begin{align}
B_i^{k + 1} = B_i^k - \frac{{B_i^k{s^k}{{({s^k})}^T}B_i^k}}{{{{({s^k})}^T}B_i^k{s^k}}} + \frac{{y_i^k{{(y_i^k)}^T}}}{{{{({s^k})}^T}y_i^k}},
\end{align}

where ${s^k} = {x^{k + 1}} - {x^k} = {\lambda _k}{d^k}$, $y_i^k = \nabla {g_i}({x^{k + 1}}) - \nabla {g_i}({x^k})$.

Step 5 Set $k: = k + 1$, and go to Step 2.\\

\textbf{Remark 3.6}~~(1) It follows  from Lemma 3.4 that  Algorithm 3.1 stops at Step 3 with a Pareto stationary point or produces an infinite sequence of nonstationary points $\{ {x^k}\} $. If Step 4 is reached in some iteration $k$, it means that in Step 3, ${d^k} \ne 0$, or equivalently, $ {\alpha}({x^k}) < 0.$   It follows from the Armijo condition  and Lemma 3.6 that objective value  sequence $\{{F}({x^k}) \}$ is $\mathbb{R}_{++}^m$-decrease, i.e., \[{F}({x^{k+1}} ) < {F}({x^k})  ~{\rm for ~all} ~k.\]

(2) It follows from \cite{MBP} and the references therein that if $g_i$ is a strongly convex function, then the matrix ${B_i^{k + 1}}$ obtained from each of the mentioned updating formula (21) for approximating the Hessian matrix $\nabla^2 g_i(x^{k + 1})$   always preserves positive definiteness.

\vskip 0.4in

\section{  Convergence analysis                }
\vskip 0.2in

In this section, we prove that the sequences generated by  Algorithm  3.1    converges to Pareto stationary points    of (CMOP) and the  Condition 3.1  will be replaced by   the following form which has been used in  \cite{MS}:

{\bf Condition 4.1}     For all $k$ and all $j=1,...,m$, we have $z^T B_j(x^k)z\geq \sigma \|z\|^2$.
\vskip 0.4in


\textbf{Lemma 4.1}~Suppose that Condition 4.1 holds true,  $\{ {d^k}\} $ is generated by Algorithm 3.1 and   that $\{ {f_i}({x^k})\} $ is bounded from below for all $i = 1,...,m$. Then, it follows that
\[\mathop {\lim }\limits_{k \to \infty } {\lambda _k}{\left\| {{d^k}} \right\|^2} = 0.\]

\begin{proof}[\quad\quad\textbf{Proof}] It follows from Lemma 3.6 and step 4 of Algorithm 3.1   that there exists $\lambda_k>0$ such that  for each $i \in \Xi_m$,

\[{f_i}({x^k} + {\lambda _k}{d^k}) \le {f_i}({x^k}) + \lambda_k \rho \alpha(x^k).\]

Adding up the above inequality from $k = 0$ to $k = \hat k$, where $\hat k$ is a positive integer, we get
\begin{align}
~{f_i}({x^{\hat k + 1}}) \le {f_i}({x^0}) + \rho \sum\limits_{k = 0}^{\hat k} {\lambda _k}  \alpha(x^k).
\end{align}

By (22) and Theorem 3.1, we get

\begin{align}
~{f_i}({x^{\hat k + 1}}) \le {f_i}({x^0}) -\frac{\sigma}{2} \rho \sum\limits_{k = 0}^{\hat k} {\lambda _k}  {\left\| {{d^k}} \right\|^2}.
\end{align}

Since $\{ {f_i}({x^k})\} $ is bounded from below for all $i = 1,...,m$,  there exists ${\hat f_i} \in \mathbb{R}$ such that ${\hat f_i} \le {f_i}({x^k})$ for all $i$ and $k$.

It follows from (23) that
\begin{align}
\sum\limits_{k = 0}^{\hat k} {{\lambda _k}{{\left\| {{d^k}} \right\|}^2}}  \le & \frac{2}{{\rho\sigma }}({f_i}({x^0}) - {f_i}({x^{\hat k + 1}}))\notag \\
\le & \frac{2}{{\rho \sigma}}({f_i}({x^0}) - {\hat f_i}). \notag
\end{align}

Taking $\hat k \to \infty $, we have $\sum\limits_{k = 0}^\infty  {{\lambda _k}{{\left\| {{d^k}} \right\|}^2}}  < \infty $

and hence $\mathop {\lim }\limits_{k \to \infty } {\lambda _k}{\left\| {{d^k}} \right\|^2} = 0$.
\end{proof}

\textbf{Theorem 4.1} (i)  Assume that the  Condition 4.1    holds true and ~ that $\{ {f_i}({x^k})\} $ is bounded from below for all $i = 1,...,m$, then every accumulation point of the sequence $\{ {x^k}\} $ generated by Algorithm 3.1, if it exists, is both a Pareto stationary point and a Pareto optimum of (CMOP).

 (ii) Moreover, if the level set of ${F}$ in the sense that $\{ x \in {\mathbb{R}^n}\mid F(x) \le F({x^0})\} $ is bounded, then $\{ {x^k}\} $ has accumulation points and they are all Pareto stationary points and Pareto optimums of (CMOP).

\begin{proof}[\quad\quad\textbf{Proof}] (i)  Let $\bar x$ be an accumulation point of $\{ {x^k}\} $ and let $\{ {x^{{k_j}}}{\rm{\} }}$ be a subsequence converging to $\bar x$.  We now prove that ${d}(\bar x) = 0$. On the contrary, we assume that ${d}(\bar x) \ne 0$. Then, by Lemma 4.1, we get   ${\lambda _{{k_j}}} \to 0 $. It follows from the definition of ${\lambda _{{k_j}}}$ in Step 4 of Algorithm 3.1  that for sufficiently large $j$ there exists some ${i_{{k_j}}} \in \{ 1,...,m\}$ such that
\[{f_{{i_{{k_j}}}}}({x^{{k_j}}} + {\zeta ^{ - 1}}{\lambda _{{k_j}}}{d^{{k_j}}}) > {f_{{i_{{k_j}}}}}({x^{{k_j}}}) + {\zeta ^{ - 1}}{\lambda _{{k_j}}}\rho  \alpha(x^{k_j}).\]

Since $i$ only takes finite number of values in $\{ 1,...,m\} $, we can assume that ${i_{{k_j}}} = \bar i $ without loss of generality. We thus obtain
\begin{align}
\frac{{{f_{\bar i}}({x^{{k_j}}} + {\zeta ^{ - 1}}{\lambda _{{k_j}}}{d^{{k_j}}}) - {f_{\bar i}}({x^{{k_j}}})}}{{{\zeta ^{ - 1}}{\lambda _{{k_j}}}}} > \rho   \alpha(x^{k_j}).
\end{align}

It follows from $0 < {\zeta ^{ - 1}}{\lambda _{{k_j}}} < 1$, the definition   of ${\alpha(x)}$ and Theorem 23.1 in  \cite{R} that
\begin{align}
  \alpha(x^{k_j}) \ge & \nabla {g_{\bar i}}{({x^{{k_j}}})^T}{d^{{k_j}}} + \frac{1}{2}{({d^{{k_j}}})^T}{B_{\bar i}}({x^{{k_j}}})({d^{{k_j}}}) + {h_{\bar i}}({x^{{k_j}}} + {d^{{k_j}}}) - {h_{\bar i}}({x^{{k_j}}})  \notag \\
\ge & \frac{{{\zeta ^{ - 1}}{\lambda _{{k_j}}}\nabla {g_{\bar i}}{{({x^{{k_j}}})}^T}{d^{{k_j}}} + \frac{1}{2}{\zeta ^{ - 1}}{\lambda _{{k_j}}}{{({d^{{k_j}}})}^T}{B_{\bar i}}({x^{{k_j}}})({d^{{k_j}}}) + {h_{\bar i}}({x^{{k_j}}} + {\zeta ^{ - 1}}{\lambda _{{k_j}}}{d^{{k_j}}}) - {h_{\bar i}}({x^{{k_j}}})}}{{{\zeta ^{ - 1}}{\lambda _{{k_j}}}}}  \notag \\
= & \frac{{{g_{\bar i}}({x^{{k_j}}} + {\zeta ^{ - 1}}{\lambda _{{k_j}}}{d^{{k_j}}}) + {h_{\bar i}}({x^{{k_j}}} + {\zeta ^{ - 1}}{\lambda _{{k_j}}}{d^{{k_j}}}) - {g_{\bar i}}({x^{{k_j}}}) - {h_{\bar i}}({x^{{k_j}}}) + o({{({\zeta ^{ - 1}}{\lambda _{{k_j}}}\left\| {{d^{{k_j}}}} \right\|)}^2})}}{{{\zeta ^{ - 1}}{\lambda _{{k_j}}}}}  \notag \\
= & \frac{{{f_{\bar i}}({x^{{k_j}}} + {\zeta ^{ - 1}}{\lambda _{{k_j}}}{d^{{k_j}}}) - {f_{\bar i}}({x^{{k_j}}})}}{{{\zeta ^{ - 1}}{\lambda _{{k_j}}}}} + \frac{{o({{({\zeta ^{ - 1}}{\lambda _{{k_j}}}\left\| {{d^{{k_j}}}} \right\|)}^2})}}{{{\zeta ^{ - 1}}{\lambda _{{k_j}}}}}, \notag
\end{align}
where ${B_{\bar i}}({x^{{k_j}}})$ is an approximation of ${\nabla ^2}{g_{\bar i}}({x^{{k_j}}})$  which updated by (21). Thus, we obtain
\begin{align}
  \alpha(x^{k_j}) \ge \frac{{{f_{\bar i}}({x^{{k_j}}} + {\zeta ^{ - 1}}{\lambda _{{k_j}}}{d^{{k_j}}}) - {f_{\bar i}}({x^{{k_j}}})}}{{{\zeta ^{ - 1}}{\lambda _{{k_j}}}}} + \frac{{o({{({\zeta ^{ - 1}}{\lambda _{{k_j}}}\left\| {{d^{{k_j}}}} \right\|)}^2})}}{{{\zeta ^{ - 1}}{\lambda _{{k_j}}}}}.
\end{align}

It follows from (24) and (25) that
\[\frac{{{f_{\bar i}}({x^{{k_j}}} + {\zeta ^{ - 1}}{\lambda _{{k_j}}}{d^{{k_j}}}) - {f_{\bar i}}({x^{{k_j}}})}}{{{\zeta ^{ - 1}}{\lambda _{{k_j}}}}} > \rho \frac{{{f_{\bar i}}({x^{{k_j}}} + {\zeta ^{ - 1}}{\lambda _{{k_j}}}{d^{{k_j}}}) - {f_{\bar i}}({x^{{k_j}}})}}{{{\zeta ^{ - 1}}{\lambda _{{k_j}}}}} + \rho \frac{{o({{({\zeta ^{ - 1}}{\lambda _{{k_j}}}\left\| {{d^{{k_j}}}} \right\|)}^2})}}{{{\zeta ^{ - 1}}{\lambda _{{k_j}}}}}.\]

Therefore, we have
\begin{align}
\frac{{{f_{\bar i}}({x^{{k_j}}} + {\zeta ^{ - 1}}{\lambda _{{k_j}}}{d^{{k_j}}}) - {f_{\bar i}}({x^{{k_j}}})}}{{{\zeta ^{ - 1}}{\lambda _{{k_j}}}}} > (\frac{\rho }{{1 - \rho }})\frac{{o({{({\zeta ^{ - 1}}{\lambda _{{k_j}}}\left\| {{d^{{k_j}}}} \right\|)}^2})}}{{{\zeta ^{ - 1}}{\lambda _{{k_j}}}}}.
\end{align}

It follows from Theorem 3.1  that
\[  \alpha(x^{k_j}) \le  - \frac{\sigma}{2} {\left\| {{d^{{k_j}}}} \right\|^2}.\]

Since ${d^{{k_j}}} \to {d}(\bar x) \ne 0$, by the above inequality, (25) and (26), we know that  there exists   $\gamma = \frac{\sigma}{2}  {\left\| d(\bar x) \right\|^2}   > 0$ such that

$$
- \gamma  \ge   \lim_{j\rightarrow \infty}   \alpha(x^{k_j}) \notag \\
\ge    \lim_{j\rightarrow \infty}  [\frac{{{f_{\bar i}}({x^{{k_j}}} + {\zeta ^{ - 1}}{\lambda _{{k_j}}}{d^{{k_j}}}) - {f_{\bar i}}({x^{{k_j}}})}}{{{\zeta ^{ - 1}}{\lambda _{{k_j}}}}} + \frac{{o({{({\zeta ^{ - 1}}{\lambda _{{k_j}}}\left\| {{d^{{k_j}}}} \right\|)}^2})}}{{{\zeta ^{ - 1}}{\lambda _{{k_j}}}}}]
$$
$$\ge    \lim_{j\rightarrow \infty}  [(\frac{\rho }{{1 - \rho }})\frac{{o({{({\zeta ^{ - 1}}{\lambda _{{k_j}}}\left\| {{d^{{k_j}}}} \right\|)}^2})}}{{{\zeta ^{ - 1}}{\lambda _{{k_j}}}}} + \frac{{o({{({\zeta ^{ - 1}}{\lambda _{{k_j}}}\left\| {{d^{{k_j}}}} \right\|)}^2})}}{{{\zeta ^{ - 1}}{\lambda _{{k_j}}}}}] =0,$$

   which contradicts the fact that $\gamma  > 0$. Therefore, we conclude that ${d}(\bar x) = 0$.

It follows from   Lemma 3.4 that $\bar x$ is Pareto stationary point of (CMOP).  By Lemma 2.1 and the strong convexity of $f_i$, we know that $\bar x$ is also a Pareto optimum of (CMOP).

(ii)It follows from  Remark 4.1(1)  that the objective value sequence $\{{F}({x^k}) \}$ is $\mathbb{R}_{++}^m$-decrease. Moreover,  since the set  $\{ x \in {\mathbb{R}^n}\mid F(x) \le F({x^0})\} $ is bounded, Thus, the sequence $\{ {x^k}\} $ generated by Algorithm 3.1 is contained in the above set and so it is also bounded and has at least one accumulation point, which is a Pareto stationary point and a Pareto optimum of (CMOP) according to the first statement.
\end{proof}

\vskip 0.4in

\section{     Conclusion      }
\vskip 0.2in

First, for the composite multiobjective optimization problems (in short, CMOP), where each objective function is the sum of a twice continuously differentiable strongly convex function and a proper convex and lower semicontinuous  but not necessarily differentiable function, the BFGS  quasi-Newnon method   with Armijo line search       are introduced. Secondly, under appropriate conditions, we prove that each cluster point of the sequence generated by the BFGS  quasi-Newton   algorithms is both a Pareto stationary point and a Pareto optimum of (CMOP). Thirdly,    numerical experiments are performed to verify the effectiveness of the proposed algorithms. In the future, there still exists an interesting topics about  the convergence rate of the proposed algorithms.

\renewcommand\refname


{\flushleft References}
\begin{thebibliography}{99}
\makeatletter
\renewcommand\@biblabel[1]{#1.}
\makeatother


\bibitem{Ansary} M. A. T. Ansary, A Newton-type proximal gradient method for nonlinear
multi-objective optimization problems, OPTIMIZATION METHODS $\&$ SOFTWARE
https://doi.org/10.1080/10556788.2022.2157000


\bibitem{AE}  J. P. Aubin, and I. Ekeland, {\em Applied Nonlinear
Analisis}, John Wiley $\&$Sons (1984).

\bibitem{Bagirov}     A. Bagirov,  N. Karmitsa
M. M. M\"{a}kel\"{a},  Introduction to Nonsmooth Optimization, Theory, Practice and Software, Springer International Publishing, Switzerland,  2014

\bibitem{BLM} Bello Cruz, J.Y., Lucambio P¨¦rez, L.R., Melo, J.G.: Convergence of the projected gradient method for quasiconvex multiobjective optimization. Nonlinear Anal. {\bf 74}(16), 5268--5273 (2011)

\bibitem{B1} Bertsekas, D.P.: Nonlinear Programming, 2nd edn. Athena Scientific, Belmont (1999)
\bibitem{BIS} Bonnel, H., Iusem, A.N., Svaiter, B.F.: Proximal methods in vector optimization. SIAM J. Optim. {\bf 15}(4), 953--970 (2005)

\bibitem{CT} Chen, G., Teboulle, M.: Convergence analysis of a proximal-like minimization algorithm using Bregman functions. SIAM J. Optim. {\bf 3}(3), 538--543 (1993)
\bibitem{DDFL} Da Cruz Neto, J.X., Da Silva, G.J.P., Ferreira, O.P., Lopes, J.O.: A subgradient method for multiobjective optimization. Comput. Optim. Appl. {\bf 54}(3), 461--472 (2013)
\bibitem{EG}  Eichfelder G.: Adaptive Scalarization Methods in Multiobjective Optimization[M], Springer -Verlag Berlin Heidelberg, 2008
\bibitem{EM}  Ehrgott, M.: Multicriteria Optimization[M],  Springer, Berlin Heidelberg, 2005

\bibitem{FGS} Fliege, J., Grana Drummond, L.M., Svaiter, B.F.: Newton¡¯s method for multiobjective optimization. SIAM J. Optim. {\bf 20}(2), 602--626 (2009)
\bibitem{FS} Fliege, J., Svaiter, B.F.: Steepest descent methods for multicriteria optimization. Math. Methods Oper. Res. {\bf 51}(3), 479--494 (2000)
\bibitem{FG} Fukuda, E.H., Grana Drummond, L.M.: A survey on multiobjective descent methods. Pesquisa Operacional {\bf 34}(3), 585--620 (2014)

\bibitem{H} Hogan, W.W.: Point-to-set maps in mathematical programming. SIAM Rev. {\bf 15}(3), 591--603 (1973)
\bibitem{MBP}   Morovati   V.,   Basirzadeh H.,
 Pourkarimi L.: Quasi-Newton methods for multiobjective optimization
problems, 4OR,   16(3), 261¨C294(2018)
\bibitem{NW} Nocedal, J., Wright, S.J.: Numerical Optimization, 2nd ed. Springer Science and Business Media, LLC, New York (2006)
\bibitem{NYY} Nakayama H., Yun  Y., Yoon M.: Sequential Approximate Multiobjective Optimization Using Computational Intelligence[M]. Springer -Verlag Berlin Heidelberg, 2009
 \bibitem{PRY}   Peng J. W., Ren J.£¬¡¡J. C.¡¡Yao,     Proximal Quasi-Newton Methods for the Composite Multiobjective Optimization Problems,  J. Nonlinear Convex Anal., To appear


\bibitem{P} Povalej, Z.: Quasi-Newton's method for multiobjective optimization. J. Comput. Appl. Math. {\bf 255}, 765--777 (2014)

\bibitem{QGC} S. Qu, M. Goh, and F.T.S. Chan, Quasi-Newton methods for solving multiobjective optimization,
Oper. Res. Lett. {bf  39},  397-¨C399 (2011)
\bibitem{R} Rockafellar, R.T.: Convex Analysis. University Press, Princeton (1970)

\bibitem{TFY} Tanabe, H., Fukuda, E.H., Yamashita, N.: Proximal gradient methods for multiobjective optimization and their applications. Comput. Optim. Appl. {\bf 72}, 339--361 (2019)

\bibitem{RP}     Ren, J., Peng, J.W.: Proximal Newton Methods for Multiobjective Optimization Problems. Acta Mathematicae Applicatae Sinica. {\bf 45}(2), 222--237 (2022)


\bibitem{MS}
 Mahdavi-Amiri, N.,   Salehi Sadaghiani, F.: A superlinearly convergent nonmonotone quasi-Newton
method for unconstrained multiobjective optimization, Optimization  Methods  $\&$ Software, {\bf 35}(6) 1223-1247 (2020)
\end{thebibliography}
\end{document}